\documentclass[a4paper,10pt,final]{article}

\usepackage{amsthm,amsfonts,amsmath}
\usepackage{mathrsfs}
\usepackage{color}
\usepackage{hyperref}
\usepackage{a4wide}
\usepackage{enumerate,float}

\usepackage{pgf,tikz}
\usetikzlibrary{arrows}
\usetikzlibrary{patterns}

\newtheorem{proposition}{Proposition}[section]
\newtheorem{theorem}[proposition]{Theorem}

\newenvironment{proofof}[1]{\smallskip\noindent{\textbf{Proof~of~#1.}}%
  \hspace{1pt}}{\hspace{-5pt}{\nobreak\quad\nobreak\hfill\nobreak%
    $\square$\vspace{2pt}\par}\smallskip\goodbreak}

\numberwithin{equation}{section}
\allowdisplaybreaks[4]

\renewcommand{\phi}{\varphi}
\renewcommand{\epsilon}{\varepsilon}
\renewcommand{\theta}{\vartheta}
\renewcommand{\L}[1]{\mathbf{L^#1}}
\newcommand{\Lloc}[1]{\mathbf{L^{#1}_{loc}}}
\newcommand{\C}[1]{\mathbf{C^{#1}}}

\newcommand{\BV}{\mathbf{BV}}
\newcommand{\modulo}[1]{{\left|#1\right|}}
\newcommand{\norma}[1]{{\left\|#1\right\|}}
\newcommand{\reali}{{\mathbb{R}}}

\newcommand{\tv}{\mathop\mathrm{TV}}
\renewcommand{\O}{\mathinner{\mathcal{O}(1)}}

\newcommand{\pint}[1]{\mathaccent23{#1}}

\renewcommand{\d}[1]{\mathinner{\mathrm{d}{#1}}}

\renewcommand{\hat}[1]{\widehat{#1}}

\newcommand{\LC}{\mathbf{LC}}

\begin{document}

\title{Uniqueness of the $1$D Compressible to Incompressible Limit}

\author{Rinaldo M. Colombo$^1$ \and Graziano Guerra$^2$}

\footnotetext[1]{INDAM Unit, University of Brescia,
  Italy. \texttt{rinaldo.colombo@unibs.it}}

\footnotetext[2]{Department of Mathematics and Applications,
  University of Milano - Bicocca,
  Italy. \texttt{graziano.guerra@unimib.it}}

\maketitle

\begin{abstract}
  \noindent Consider two compressible immiscible fluids in 1D in the
  isentropic approximation. The first fluid is surrounded and in
  contact with the second one. As the Mach number of the first fluid
  vanishes, the coupled dynamics of the two fluids results as the
  compressible to incompressible limit and is known to satisfy an
  ODE--PDE system. Below, a characterization of this limit is
  provided, ensuring its uniqueness.

  \medskip

  \noindent\textbf{Keywords:} Compressible to Incompressible limit,
  Hyperbolic Conservation Laws, Uniqueness of the Zero Mach number
  Limit

  \medskip

  \noindent\textbf{2010 MSC:} 35L65, 35Q35, 76N99
\end{abstract}

\section{Introduction}
\label{sec:Intro}

The literature on the compressible to incompressible limit is vast. We
refer for instance to the well known
results~\cite{KlainermanMajda1981, KlainermanMajda1982,
  MetivierSchochet2001, Schochet1986}, the more
recent~\cite{ChenCleopatra2008, JiangYong}, the
review~\cite{Schochet2005} and the references therein.

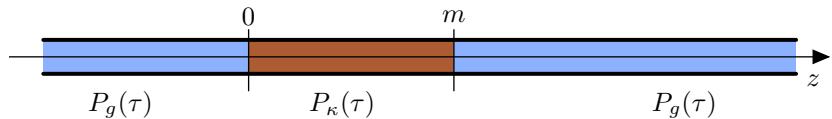
\begin{figure}[!ht]
  \centering \definecolor{wwzzff}{rgb}{0.4,0.6,1.0}
  \definecolor{zzttqq}{rgb}{0.6,0.2,0.0}
  \definecolor{cqcqcq}{rgb}{0.7529411764705882,0.7529411764705882,0.7529411764705882}
  \begin{tikzpicture}[line cap=round,line join=round,>=triangle
    45,x=0.45cm,y=0.45cm]
    \clip(-4.00000000000001,-4.00000000000001) rectangle
    (20.0000000000003,-0.0); \fill[color=zzttqq,fill=zzttqq,fill
    opacity=0.8] (3.0,-1.0) -- (3.0,-2.0) -- (9.0,-2.0) -- (9.0,-1.0)
    -- cycle; \fill[color=wwzzff,fill=wwzzff,fill opacity=0.75]
    (-3.0,-1.0) -- (-3.0,-2.0) -- (3.0,-2.0) -- (3.0,-1.0) -- cycle;
    \fill[color=wwzzff,fill=wwzzff,fill opacity=0.75] (9.0,-1.0) --
    (9.0,-2.0) -- (19.0,-2.0) -- (19.0,-1.0) -- cycle; \draw [line
    width=1.5000000000000003pt] (-3.0,-1.0)-- (19.0,-1.0); \draw [line
    width=1.5000000000000003pt] (19.0,-2.0)-- (-3.0,-2.0);
    \draw[color=black] (-0.74,-2.900000000000002) node
    {$P_{g}(\tau)$}; \draw[color=black] (15.74,-2.900000000000002)
    node {$P_{g}(\tau)$}; \draw[color=black] (5.74,-2.900000000000002)
    node {$P_{\kappa}(\tau)$}; \draw[color=black]
    (3.0,-0.3000000000002) node {$0$}; \draw[color=black]
    (9.0,-0.3000000000002) node {$m$}; \draw [line
    width=0.5000000000000003pt,<-] (20.0,-1.5)-- (-4.0,-1.5);
    \draw[color=black] (19.5,-2.2000000000002) node {$z$}; \draw [line
    width=0.5000000000000003pt] (3.00,-2.5)-- (3.00,-0.7); \draw [line
    width=0.5000000000000003pt] (9.00,-2.5)-- (9.00,-0.7);
  \end{tikzpicture}
  \caption{In Lagrangian coordinates, the boundaries separating the
    two fluids are fixed.}
  \label{fig:liquidogas2}
\end{figure}

In this paper, following~\cite{ColomboGuerraIncompressible}, we
consider two compressible immiscible fluids and study the limit as one
of the two becomes incompressible. A volume of a compressible inviscid
fluid, say the \emph{liquid}, is surrounded by another compressible
fluid, say the \emph{gas}. Using the Lagrangian formulation, in the
isentropic case, we assume that the gas obeys a fixed pressure law
$P_{g}(\tau)$, while for the liquid we assume a one parameter family
of pressure laws $P_{\kappa}(\tau)$ such that
$P_{\kappa}^{\prime}(\tau) \to -\infty$ as $\kappa\to 0$. The total
mass of the liquid is fixed so that in Lagran\-gian coordinates the
liquid and gas phases fill the fixed sets (see
Figure~\ref{fig:liquidogas2})
\begin{displaymath}
  \mathcal{L} = \left]0,m \right[
  \qquad \mbox{ and } \qquad
  \mathcal{G} = \reali \setminus \left]0,m \right[ \,.
\end{displaymath}
For an Eulerian description, see~\cite{ColomboGuerraIncompressible}.

On $P_{g}\left(\tau\right)$ and $P_{\kappa}\left(\tau\right)$, we
require the usual hypotheses and the incompressible limit assumption:
\begin{equation}
  \label{eq:Pproperties}
  P_g, P_\kappa \in \C4,\,
  P_{g}\left(\tau\right),\, P_{\kappa}\left(\tau\right)>0;\;
  P'_{g}\left(\tau\right),\, P_{\kappa}^{\prime}\left(\tau\right)<0;\;
  P''_{g}\left(\tau\right),\, P_{\kappa}^{\prime\prime}\left(\tau\right)>0;\;
  P_{\kappa}^{\prime}\left(\tau\right)\xrightarrow{\kappa\to 0}-\infty \,.
\end{equation}
The standard choice $P_g (\tau) = k / \tau^\gamma$
satisfies~\eqref{eq:Pproperties} for all $k > 0$ and $\gamma >0$.

The coupled dynamics of the two fluids is described by the
$p$-system~\cite[Formula~(7.1.11)]{DafermosBook}
\begin{equation}
  \label{eq:baseeqint}
  \begin{cases}
    \partial_{t}\tau-\partial_{z}v=0
    \\
    \partial_{t}v+\partial_{z}P_{\kappa}\left(z,\tau\right)=0,
  \end{cases}
  \quad \mbox{ where } \quad
  P_{\kappa}\left(z,\tau\right)=
  \begin{cases}
    P_{\kappa}\left(\tau\right) & \mbox{ for }z\in \mathcal{L}
    \\
    P_{g}\left(\tau\right) & \mbox{ for }z\in \mathcal{G},
  \end{cases}
\end{equation}
$v(t,z)$ being the fluid speed at time $t$ and at the Lagrangian
coordinate $z$.

In Lagrangian coordinates, the conservation of mass and momentum are
equivalent to the conservation of $\tau$ and $v$ which, in turn, are
equivalent along the interfaces $z=0$ and $z=m$ to the
Rankine--Hugoniot conditions for~\eqref{eq:baseeqint}. Therefore, for
a.e.~$t\geq 0$,
\begin{displaymath}
  \begin{cases}
    v\left(t,0-\right)=v\left(t,0+\right)
    \\
    P_{g}\left(\tau\left(t,0-\right)\right) =
    P_{\kappa}\left(\tau\left(t,0+\right)\right),
  \end{cases}
  \quad
  \begin{cases}
    v\left(t,m-\right)=v\left(t,m+\right)
    \\
    P_{\kappa}\left(\tau\left(t,m-\right)\right) =
    P_{g}\left(\tau\left(t,m+\right)\right).
  \end{cases}
\end{displaymath}
In other words, pressure and velocity have to be continuous across the
interfaces. Hence, the pressure is a natural choice as unknown, rather
than the specific volume.
Following~\cite{ColomboGuerraIncompressible, ColomboGuerraInc,
  ColomboGuerraSchleper, CGSInc}, we introduce the inverse functions
of the pressure laws
\begin{equation}
  \label{eq:Tproperties}
  \mathcal{T}_{g}(p)=P^{-1}_{g}\left(p\right),\quad
  \mathcal{T}_{\kappa}(p)=P^{-1}_{\kappa}\left(p\right)
  \quad \mbox{ where} \quad
  \mathcal{T}_{\kappa}'\left(p\right)\xrightarrow{\kappa\to 0} 0 \,,
\end{equation}
the last limit being a consequence of~\eqref{eq:Pproperties}. Rewrite
system~\eqref{eq:baseeqint} with $(p,v)$ as unknowns
\begin{equation}
  \label{eq:baseeqinp}
  \begin{cases}
    \partial_{t}\mathcal{T}_{\kappa}\left(z,p\right)-\partial_{z}v=0
    \\
    \partial_{t}v+\partial_{z}p=0 \,,
  \end{cases}
  \quad \mbox{ where } \quad
  \mathcal{T}_{\kappa}\left(z,p\right)=
  \begin{cases}
    \mathcal{T}_{\kappa}\left(p\right) & \mbox{ for }z\in \mathcal{L}
    \\
    \mathcal{T}_{g}\left(p\right) & \mbox{ for }z\in\mathcal{G} \,.
  \end{cases}
\end{equation}
The conditions at the interfaces become continuity requirements on the
unknown functions:
\begin{equation}
  \label{eq:intCond}
  \begin{cases}
    v\left(t,0-\right)=v\left(t,0+\right)
    \\
    p\left(t,0-\right)=p\left(t,0+\right)
  \end{cases}
  \qquad
  \begin{cases}
    v\left(t,m-\right)=v\left(t,m+\right)
    \\
    p\left(t,m-\right)=p\left(t,m+\right)
  \end{cases}\quad \mbox{for a.e. }t\geq 0 \,.
\end{equation}

As in~\cite{ColomboGuerraIncompressible}, we fix a pressure law $P$
and choose $\mathcal{T}=P^{-1}$, so that
\begin{equation}
  \label{eq:tkfamily}
  \mathcal{T}_{\kappa}\left(p\right) =
  \mathcal{T}\left(\bar p + \kappa^{2}
    \left(p-\bar p\right)\right)\,,\quad
  \lim_{\kappa\to
    0}\mathcal{T}_{\kappa}\left(p\right)=\mathcal{T}\left(\bar
    p\right) = \bar \tau\,,
\end{equation}
where $\bar \tau$ is the constant specific volume at the
incompressible limit and $\bar p = P (\bar\tau)$.  For instance, the
(modified) Tait equation of state~\cite[Formula~(1)]{MacDonald} fits
into~\eqref{eq:tkfamily} with
\begin{displaymath}
  \mathcal{T} (p) = p^{-1/n}
  \quad \mbox{ with }\quad
  \kappa^2 = \frac{n \, \beta_o}{\bar\tau^n}
\end{displaymath}
where $\beta_o$ is the isothermal compressibility, $n$ is a pressure
independent parameter and $\beta_o \to 0$ at the incompressible limit.

The main result in~\cite{ColomboGuerraIncompressible} states the
rigorous convergence (up to a subsequence) at the incompressible limit
in the liquid phase of the solutions to~\eqref{eq:baseeqinp} to
solutions to
\begin{equation}
  \label{eq:IncompL}
  \left\{
    \begin{array}{lr@{\,}c@{\,}l@{\qquad\qquad}l@{}}
      \left\{
      \begin{array}{@{}l}
        \partial_t \mathcal{T}_{g}(p) - \partial_z v = 0
        \\
        \partial_t v
        +
        \partial_z  p
        =
        0
      \end{array}
      \right.
   & z
   & \in
   &  \mathcal{G}
   & \mbox{gas;}
      \\[20pt]
      \dot v_\ell
      =
      \frac{p(t, 0-) - p(t, m+)}{m}
   & & &
   & \mbox{liquid;}
      \\[10pt]
      \left\{
      \begin{array}{@{}rcl@{}}
        v \left(t, 0-\right)
        & =
        & v_\ell (t)
        \\
        v \left(t, m+\right)
        & =
        & v_\ell (t)
      \end{array}
          \right.
   & & &
     & \mbox{interface.}
    \end{array}
  \right.
\end{equation}
The existence of a Lipschitz continuous semigroup generated
by~\eqref{eq:IncompL} is proved in~\cite{BorscheColomboGaravello2}. On
the other hand, a characterization yielding the uniqueness of
solutions to~\eqref{eq:IncompL} is obtained
in~\cite{ColomboGuerraAML}.

In this paper we show that the incompressible limit obtained
in~\cite{ColomboGuerraIncompressible} satisfies the characterization
in~\cite{ColomboGuerraAML}. Hence, the solution $(p_\kappa, v_\kappa)$
to~\eqref{eq:baseeqinp} converges as $\kappa \to 0$, the limit being
the unique solution to~\eqref{eq:IncompL}.

\medskip

The next Section is devoted to the formal statements, while
Section~\ref{sec:TP} contains the technical proofs.

\section{Main Result}
\label{sec:Main}

Throughout, we denote by $\LC$ the set of functions defined on
$\reali \setminus \left]0,m\right[$ that are locally constant out of a
compact set, i.e., they attain a constant value on
$\left]-\infty, -M\right]$ and a, possibly different, constant value
on $\left[M, +\infty\right[$, for a suitable positive $M$.

Below, solutions to~\eqref{eq:IncompL} are understood in the sense
of~\cite[Definition~3.2]{ColomboGuerraIncompressible}, see
also~\cite[Definition~2.5]{BorscheColomboGaravello2}, and are
constructed in~\cite{ColomboGuerraIncompressible} as limits of
solutions to~\eqref{eq:baseeqint}. In solutions
to~\eqref{eq:baseeqint}, the propagation speed of waves in the gas
region $\mathcal{G}$ is uniformly bounded, independently of
$\kappa$. Therefore, to prove the uniqueness of solutions
to~\eqref{eq:IncompL} obtained as the compressible to incompressible
limit, it is sufficient to consider initial data
$\left((\tau_o,v_o),v_{\ell,o}\right)$ such that $(\tau_o,v_o)$ is in
$\LC$ and $v_{\ell,o} \in \reali$.

Given
$\left((\tau,v),v_\ell\right) \in \BV (\reali;\reali^2) \times \reali$
such that $(\tau,v) \in \LC$, call
\begin{displaymath}
  (\tau_{\pm\infty}, v_{\pm\infty}) = \lim_{x\to\pm\infty} (\tau,v) (x) \,.
\end{displaymath}
Under the transformation
\begin{equation}
  \label{eq:6}
  U (x)
  =
  \left[
    \begin{array}{c}
      \tau (-x) - \tau_{-\infty}
      \\
      v (-x) - v_{-\infty}
      \\
      \tau (x+m) - \tau_{+\infty}
      \\
      v (x+m) - v_{+\infty}
    \end{array}
  \right]
  \qquad
  w
  =
  \left[
    \begin{array}{c}
      v_\ell - v_{-\infty}
      \\
      v_\ell - v_{+\infty}
    \end{array}
  \right] \,,
\end{equation}
setting
\begin{equation}
  \label{eq:7}
  \begin{array}{@{}rcl@{\quad}rcl@{}}
    f(U)
    & =
    & \left[
      \begin{array}{c}
        U_2 \\ -P_g (U_1+\tau_{-\infty}) \\ -U_4 \\ P_g (U_3+\tau_{+\infty})
      \end{array}
    \right]

    &  F (U,w)
    & =
    & \dfrac{1}{m}
      \left[
      \begin{array}{c}
        P_g(U_1 + \tau_{-\infty})- P_g (U_3 + \tau_{+\infty})
        \\
        P_g(U_1 + \tau_{-\infty})- P_g (U_3 + \tau_{+\infty})
      \end{array}
    \right]
    \\

    b (U)
    & =
    & \left[
      \begin{array}{c}
        U_2 \\ U_4
      \end{array}
    \right]
    & g (w)
    & =
    & w
  \end{array}
\end{equation}
the Cauchy Problem
\begin{equation}
  \label{eq:2}
  \left\{
    \begin{array}{@{}l@{}l@{}}
      \left\{
      \begin{array}{@{}l@{}}
        \partial_t \tau - \partial_z v = 0
        \\
        \partial_t v
        +
        \partial_z  P_g (\tau)
        =
        0
      \end{array}
      \right.
   & x \in \mathcal{G}
      \\[10pt]
      \dot v_\ell
      =
      \dfrac{P_g\left(\tau(t, 0-)\right) - P_g\left(\tau(t, m+)\right)}{m}
      \\[5pt]
      \left\{
      \begin{array}{@{}rcl@{}}
        v \left(t, 0-\right)
        & =
        & v_\ell (t)
        \\
        v \left(t, m+\right)
        & =
        & v_\ell (t)
      \end{array}
          \right.
      \\
      (\tau,v)(0,x) = (\tau_o, v_o) (x)
   & x \in \mathcal{G}
      \\
      v_\ell (0) = v_{\ell,o}
    \end{array}
  \right.
\end{equation}
is formally equivalent to
\begin{equation}
  \label{eq:4}
  \left\{
    \begin{array}{ll@{}}
      \partial_t U (t,x) + \partial_x f\left(U (t,x)\right) = 0
      & x \in \reali^+
      \\
      b\left(U (t, 0+)\right) = g\left(w (t)\right)
      \\
      \dot w (t) = F\left(U (t, 0+), w (t)\right)
      \\
      U (0,x) = U_o (x)
      & x \in \reali^+
      \\
      w (0) = w_o
    \end{array}
  \right.
\end{equation}
which fits in the well posedness theory developed
in~\cite{ColomboGuerraAML}, as proved by the following Proposition.

\begin{proposition}
  \label{prop:equi}
  Let $P_g$ satisfy~\eqref{eq:Pproperties}.  Fix
  $\tau_{-\infty}, \tau_{+\infty} \in \pint{\reali}^+$.  Then,
  system~\eqref{eq:4} generates a semigroup
  $S \colon \reali^+ \times \mathcal{D} \to \mathcal{D}$ uniquely
  characterized by the properties~(i)--(iv)
  in~\cite[Theorem~4]{ColomboGuerraAML}. Moreover, for a suitable
  positive $\delta$,
  \begin{equation}
    \label{eq:tildedelta}
    \mathcal{D}
    \supseteq
    \left\{
      (U,w) \in (\L1\cap\BV) (\reali^+; \reali^4)\times \reali^2
      \colon
      \tv (U) + \norma{b\left(U (0+)\right) - g (w)} < \delta
    \right\}\,.
  \end{equation}
\end{proposition}

The above proposition leads to the main result of this paper.

\begin{theorem}
  \label{thm:unicita}
  Let $t \to \left((\tau,v),v_\ell\right) (t)$ be a solution
  to~\eqref{eq:IncompL} obtained as limit for $\kappa \to 0$ of
  solutions to~\eqref{eq:baseeqint}, with an initial datum in $\LC$
  and satisfying for all $t \in \reali_+$
  \begin{equation}
    \label{eq:8}
    \tv\left((\tau,v) (t);\mathcal{G}\right)
    +
    \norma{\left[
        \begin{array}{c}
          v (t, 0-) - v_{\ell} (t)
          \\
          v (t,m+) - v_{\ell} (t)
        \end{array}
      \right]}
    < \delta
  \end{equation}
  with $\delta$ as in~\eqref{eq:tildedelta}. Correspondingly, define
  $t \to (U,w) (t)$ as in~\eqref{eq:6}. Then,
  \begin{enumerate}
  \item for all $t \in \reali_+$, the map $t \to (U,w) (t)$ coincides
    with an orbit of the semigroup $S$ defined in
    Proposition~\ref{prop:equi}.
  \item The semigroup $S$ is defined globally in time for all initial
    data with sufficiently small total variation.
  \end{enumerate}
\end{theorem}

\noindent In the above statement, as well as below, we use the obvious
notation
\begin{displaymath}
  \tv\left((\tau,v);\mathcal{G}\right) =
  \tv\left((\tau,v);\left]-\infty,0\right]\right) +
  \tv\left((\tau,v);\left[m, +\infty\right[\right) \,.
\end{displaymath}

\section{Technical Proofs}
\label{sec:TP}

\begin{proofof}{Proposition~\ref{prop:equi}}
  On the basis of~\eqref{eq:7} and with the help
  of~\eqref{eq:Pproperties}, we verify that~\eqref{eq:4} satisfies the
  assumptions of~\cite[Theorem~4]{ColomboGuerraAML}.  With reference
  to the notation therein, set $n = 4$, $l=2$, $m=2$. Now, observe
  that~\textbf{(H1)} holds. Clearly, $f$ is of class $\C4$
  by~\eqref{eq:Pproperties}. The strict hyperbolicity
  of~\eqref{eq:IncompL} can easily be recovered through a rescaling of
  the space variable, since the different $p$--systems
  in~\eqref{eq:IncompL} interact only through the boundary,
  see~\cite[Lemma~4.1]{c-g-h-s-2009}. Besides, with standard notation,
  we have:
  \begin{displaymath}
    \begin{array}{@{}r@{\,}c@{\,}l@{}}
      \lambda_1 (U)
      & =
      & -\sqrt{-P_g' (U_1+\tau_{-\infty})}
      \\
      \lambda_2 (U)
      & =
      & -\sqrt{-P_g' (U_3+\tau_{+\infty})}
      \\
      \lambda_3 (U)
      & =
      & \sqrt{-P_g' (U_1+\tau_{-\infty})}
      \\
      \lambda_4 (U)
      & =
      & \sqrt{-P_g' (U_3+\tau_{+\infty})}
    \end{array}
    \quad
    r_1 =
    \left[
      \begin{array}{@{}c@{}}
        1 \\ \lambda_1 (U) \\0 \\0
      \end{array}
    \right]
    \
    r_2 =
    \left[
      \begin{array}{@{}c@{}}
        0 \\0 \\ 1 \\ -\lambda_2 (U)
      \end{array}
    \right]
    \
    r_3 =
    \left[
      \begin{array}{@{}c@{}}
        1 \\ \lambda_3 (U) \\ 0 \\ 0
      \end{array}
    \right]
    \
    r_4 =
    \left[
      \begin{array}{@{}c@{}}
        0 \\ 0 \\ 1 \\ -\lambda_4 (U)
      \end{array}
    \right].
  \end{displaymath}

  Concerning~\textbf{(H2)}, $b$ is clearly of class $\C4$ and
  $b (0) = 0$. Moreover,
  \begin{displaymath}
    \det \left[Db (U) \left[r_3 (U) \quad r_4 (U)\right]\right]
    =
    \det
    \left[
      \begin{array}{cc}
        \lambda_3 (U)
        & 0
        \\
        0
        & - \lambda_4 (U)
      \end{array}
    \right]
    =
    - \lambda_3 (U) \, \lambda_4 (U)
  \end{displaymath}
  and the latter expression above is non zero
  by~\eqref{eq:Pproperties}.

  Assumptions~\textbf{(H3)} and~\textbf{(H4)} are immediate
  by~\eqref{eq:7} and~\eqref{eq:Pproperties}.

  An application of~\cite[Theorem~4]{ColomboGuerraAML} yields the
  existence of a Lipschitz continuous local semigroup $\mathcal{S}$
  defined on a domain $\mathcal{D}$
  enjoying~\cite[Properties~(i)--(iv) in
  Theorem~4]{ColomboGuerraAML}. Note that~\eqref{eq:tildedelta} holds
  by~\cite[Formula~(4) and Theorem~4]{ColomboGuerraAML}.
\end{proofof}

\begin{proofof}{Theorem~\ref{thm:unicita}}
  Given $t \to \left((\tau,v),v_\ell\right) (t)$, define
  $t \to (U,w) (t)$ by means of~\eqref{eq:6}. Since
  \begin{displaymath}
    \tv \left(U (t)\right)
    +
    \norma{b\left(U (t,0+)\right) - g \left(w (t)\right)}
    =
    \tv\left((\tau,v) (t);\mathcal{G}\right)
    +
    \norma{\left[
        \begin{array}{c}
          v (t,0-) - v_{\ell} (t)
          \\
          v (t,m+) - v_{\ell} (t)
        \end{array}
      \right]}
  \end{displaymath}
  thanks to~\eqref{eq:tildedelta} we obtain that for all
  $t \in \reali_+$, $(U,w) (t) \in \mathcal{D}$, $\mathcal{D}$ being
  the domain defined in Proposition~\ref{prop:equi}.

  For $\epsilon>0$ and $\kappa>0$, call
  $(p^{\kappa,\epsilon}, v^{\kappa,\epsilon})$ the wave front tracking
  approximate solutions to~\eqref{eq:baseeqinp}, see
  also~\cite[Formula~(2.5)]{ColomboGuerraIncompressible} as defined
  in~\cite[Section~4]{ColomboGuerraIncompressible}, converging to
  $\left((P_g(\tau), v), v_\ell\right)$ first as $\epsilon \to 0$ and
  then as $\kappa\to 0$. To simplify the notation, here we omit the
  introduction of sequences and subsequences.

  In the limit $\epsilon \to 0$, by~\cite[Proof of
  Theorem~3.3]{ColomboGuerraIncompressible} we have that
  \begin{displaymath}
    \lim_{\epsilon\to 0} (p^{\kappa,\epsilon}, v^{\kappa,\epsilon}) (t)
    =
    (p^{\kappa}, v^{\kappa}) (t)
    \quad \mbox{ for all } \quad t\geq 0
    \quad \mbox{ in } \quad \Lloc1 (\reali; \reali^2) \,,
  \end{displaymath}
  where $(p^{\kappa}, v^{\kappa})$
  solves~\cite[Formula~(2.5)]{ColomboGuerraIncompressible} in the
  sense of~\cite[Definition~3.1]{ColomboGuerraIncompressible}.

  In the limit $\kappa \to 0$, we have that
  \begin{displaymath}
    \begin{array}{rclll}
      \displaystyle
      \lim_{\kappa\to 0} (p^\kappa, v^\kappa) (t, \cdot)
      & =
      & (p,v) (t, \cdot)
      & \mbox{ for all } t \geq 0
      & \mbox{in } \Lloc1 (\mathcal{G};\reali^2)
      \\
      \displaystyle
      \lim_{\kappa\to 0} v^\kappa (t, \cdot)
      & =
      & v_\ell (t)
      & \mbox{ for all } t \geq 0
      & \mbox{in } \L1 (\mathcal{L};\reali^2)
    \end{array}
  \end{displaymath}
  $v_\ell$ being independent of $z$. Introduce
  \begin{equation}
    \label{eq:ukeps}
    \begin{array}{rcl}
      \overline{v}^{\kappa,\epsilon} (t)
      & =
      & \displaystyle \frac1m \int_0^m v^{\kappa,\epsilon} (t, z) \d{z}
      \\[10pt]
      u^{\kappa, \epsilon} (t)
      & =
      & \displaystyle
        \left(
        (p^{\kappa,\epsilon}, v^{\kappa,\epsilon}) (s)_{\strut\big|\mathcal{G}}
        \,,\;
        \overline{v}^{\kappa,\epsilon} (t)
        \right)
    \end{array}
  \end{equation}
  and note that the above $\Lloc1$ convergence implies that
  \begin{equation}
    \label{eq:convergenze}
    \begin{array}{rcrclll}
      u^{\kappa,\epsilon} (t)
      & \to
      & u^\kappa (t)
      & =
      & \left(
        (p^\kappa, v^\kappa) (t)_{\strut\big|\mathcal{G}}
        \,,\;
        \displaystyle \frac1m \int_0^m v^\kappa (t, z) \d{z}
        \right)
      & \mbox{as } \epsilon\to 0
      & \mbox{for all } t\geq 0 \,,
      \\[20pt]
      u^\kappa (t)
      & \to
      & u (t)
      & =
      & \left((p, v), v_\ell\right) (t)
      &\mbox{as } \kappa\to 0
      & \mbox{for all } t\geq 0 \,.
    \end{array}
  \end{equation}
  Following~\eqref{eq:6} and~\eqref{eq:ukeps}, introduce the variables
  \begin{equation}
    \label{eq:Uw}
    U^{\kappa,\epsilon} (t,x)
    =
    \left[
      \begin{array}{c}
        \mathcal{T}_g\left( p^{\kappa,\epsilon}(t,-x) \right)- \tau_{-\infty}
        \\
        v^{\kappa,\epsilon} (t,-x) - v_{-\infty}
        \\
        \mathcal{T}_g\left(p^{\kappa,\epsilon} (t, x+m)\right) - \tau_{+\infty}
        \\
        v^{\kappa,\epsilon} (t, x+m) - v_{+\infty}
      \end{array}
    \right]
    \qquad
    w^{\kappa,\epsilon} (t)
    =
    \left[
      \begin{array}{c}
        \overline{v}^{\kappa,\epsilon} (t) - v_{-\infty}
        \\[5pt]
        \overline{v}^{\kappa,\epsilon} (t) - v_{+\infty}
      \end{array}
    \right]
  \end{equation}
  and the distance
  \begin{displaymath}
    d\left((U,w),\,(\tilde U, \tilde w)\right)
    =
    \norma{\tilde U - U}_{\L1 (\reali_+; \reali^4)}
    +
    \norma{\tilde w - w} \,.
  \end{displaymath}
  By the convergences~\eqref{eq:convergenze}, the
  definition~\eqref{eq:Uw} and the continuity of $S_t$
  \begin{equation}
    \label{eq:agognata}
    d \left((U,w) (t), S_t \left((U,w) (0)\right)\right)
    \leq
    \lim_{\kappa \to 0} \lim_{\epsilon \to 0}
    d
    \left(
      (U^{\kappa,\epsilon},w^{\kappa,\epsilon}) (t),
      S_t \left((U^{\kappa,\epsilon},w^{\kappa,\epsilon}) (0)\right)
    \right) \,.
  \end{equation}
  By~\cite[Theorem~2.9]{BressanLectureNotes}, denoting by $L$ a
  Lipschitz constant of $S_t$,
  \begin{eqnarray}
    \nonumber
    &
    &  d \left(
      (U^{\kappa,\epsilon},w^{\kappa,\epsilon}) (t),
      S_t \left((U^{\kappa,\epsilon},w^{\kappa,\epsilon}) (0)\right)
      \right)
    \\
    \nonumber
    & \leq
    & L
      \int_0^t
      \liminf_{h\to 0}
      \frac{1}{h} \;
      d\left(
      (U^{\kappa,\epsilon},w^{\kappa,\epsilon}) (s+h),
      S_h \left((U^{\kappa,\epsilon},w^{\kappa,\epsilon}) (s)\right)
      \right) \d{s}
    \\
    \label{eq:dX}
    & =
    & L
      \int_0^t
      \liminf_{h\to 0}
      \frac{1}{h} \;
      d\left(
      (U^{\kappa,\epsilon},w^{\kappa,\epsilon}) (s+h),
      \mathcal{F} (h) \left((U^{\kappa,\epsilon},w^{\kappa,\epsilon}) (s)\right)
      \right) \d{s}
  \end{eqnarray}
  where $\mathcal{F}$ is the local flow defined
  in~\cite[Formula~(5)]{ColomboGuerraAML}. By construction, the last
  term in the integrand above is
  \begin{eqnarray}
    \nonumber
    \!\!\!\!\!\!\!\!\!\!\!\!\!\!\!\!\!\!\!\!\!
    & \!\!\!
    & d\left(
      (U^{\kappa,\epsilon},w^{\kappa,\epsilon}) (s+h),
      \mathcal{F} (h) \left((U^{\kappa,\epsilon},w^{\kappa,\epsilon}) (s)\right)
      \right)
    \\
    \label{eq:questa}
    \!\!\!\!\!\!\!\!\!\!\!\!\!\!\!\!\!\!\!\!\!
    & = \!\!\!
    & \norma{U^{\kappa,\epsilon} (s+h)
      - \bar S_h \left(\bar U^{\kappa,\epsilon} (s)\right)}_{\L1 (\reali_+;\reali^4)}
      +
      \norma{
      w^{\kappa,\epsilon} (s+h)
      -
      \left[
      w^{\kappa,\epsilon} (s)
      +
      h F\left(U^{\sigma} ,w^{\kappa,\epsilon} (s)\right)\right]
      }
  \end{eqnarray}
  where $F$ is as in~\eqref{eq:7}, $\bar S$ is the Standard Riemann
  Semigroup~\cite[Chapter~9]{BressanLectureNotes} generated by
  $\partial_t U + \partial_x f (U)=0$, with $f$ as in~\eqref{eq:7},
  \begin{displaymath}
    \bar U^{\kappa,\epsilon} (s,x)
    =
    \left\{
      \begin{array}{l@{\qquad}r@{\,}c@{\,}l}
        U^{\kappa,\epsilon} (s,x)
        & x
        & \geq
        & 0
        \\
        U^\sigma
        & x
        & <
        &0
      \end{array}
    \right.
  \end{displaymath}
  and $U^\sigma$ is the unique state satisfying
  $b (U^\sigma) = g\left(w^{\kappa,\epsilon} (s)\right)$ that can be
  connected to $U^{\kappa,\epsilon} (t,0+)$ by means of Lax waves with
  positive speed, with $b$ and $g$ as in~\eqref{eq:7}.

  Introduce
  \begin{displaymath}
    (\bar p, \bar v) (z)
    =
    \left\{
      \begin{array}{l@{\qquad}r@{\,}c@{\,}l}
        (p^{\kappa,\epsilon}, v^{\kappa,\epsilon}) (s,z)
        & z
        & <
        & 0
        \\
        (p^\sigma_0, \overline{v}^{k,\epsilon}) (s)
        & z
        & \in
        & \left[0, m/2\right[
        \\
        (p^\sigma_m, \overline{v}^{k,\epsilon}) (s)
        & z
        & \in
        & [m/2, m]
        \\
        (p^{\kappa,\epsilon}, v^{\kappa,\epsilon}) (s,z)
        & z
        & >
        & m
      \end{array}
    \right.
  \end{displaymath}
  where $\overline{v}^{k,\epsilon}$ is defined in~\eqref{eq:ukeps} and
  the pressure $p^\sigma_0$, respectively $p^\sigma_m$, is such that
  the Riemann Problem
  \begin{displaymath}
    \left\{
      \begin{array}{@{}l}
        \partial_{t}\mathcal{T} (p)-\partial_{z}v = 0
        \\
        \partial_{t}v+\partial_{z}p=0
        \\
        (p,v) (0,x) =
        \left\{
        \begin{array}{@{}l@{\quad}r@{\,}c@{\,}l}
          (p^{\kappa,\epsilon}, v^{\kappa,\epsilon}) (s,0-)
          & x
          & <
          & 0
          \\
          (p^\sigma_0, \overline{v}^{\kappa,\epsilon}) (s)
          & x
          & >
          & 0 \,,
        \end{array}
            \right.
      \end{array}
    \right.
    \mbox{ resp. }
    \left\{
      \begin{array}{@{}l@{}}
        \partial_{t}\mathcal{T} (p)-\partial_{z}v = 0
        \\
        \partial_{t}v+\partial_{z}p=0
        \\
        (p,v) (0,x) =
        \left\{
        \begin{array}{@{}l@{\quad}r@{\,}c@{\,}l@{}}
          (p^\sigma_m, \overline{v}^{\kappa,\epsilon}) (s)
          & x
          & <
          & 0
          \\
          (p^{\kappa,\epsilon}, v^{\kappa,\epsilon}) (s,m+)
          & x
          & >
          & 0 \,,
        \end{array}
            \right.
      \end{array}
    \right.
  \end{displaymath}
  is solved by waves with negative, respectively positive, speed. Note
  that by~\cite[Lemma~4.1]{ColomboGuerra2010}
  \begin{equation}
    \label{eq:10}
    \begin{array}{rcl}
      \modulo{p^{\kappa,\epsilon} (s,0) - p^\sigma_0}
      & =
      & \O \,
        \modulo{v^{\kappa,\epsilon} (s,0) - \overline{v}^{\kappa,\epsilon} (s)}
        \,,
      \\[10pt]
      \modulo{p^{\kappa,\epsilon} (s,m) - p^\sigma_m}
      & =
      & \O \,
        \modulo{v^{\kappa,\epsilon} (s,m) - \overline{v}^{\kappa,\epsilon} (s)}
        \,,
    \end{array}
  \end{equation}
  recall that $z \to p^{\kappa,\epsilon} (s,z)$ and
  $z \to v^{\kappa,\epsilon} (s,z)$ are locally constant in
  neighborhoods of $z=0$ and $z=m$,
  see~\cite[Formula~(4.12)]{ColomboGuerraIncompressible}.  Call
  $\Sigma$ the Standard Riemann
  Semigroup~\cite[Chapter~9]{BressanLectureNotes} generated by the
  $p$-system $\left\{
    \begin{array}{@{}l}
      \partial_{t}\mathcal{T} (p)-\partial_{z}v = 0
      \\
      \partial_{t}v+\partial_{z}p=0
    \end{array}
  \right.$
  for $z$ varying on all the real line. Observe that the first addend
  in~\eqref{eq:questa} reads
  \begin{equation}
    \label{eq:9}
    \norma{U^{\kappa,\epsilon} (s+h)
      - \bar S_h \left(\bar U^{\kappa,\epsilon} (s)\right)}_{\L1 (\reali_+;\reali^4)}
    \leq
    \int_{\mathcal{G}}
    \norma{
      (p^{\kappa,\epsilon},v^{\kappa,\epsilon}) (s+h,z)
      -
      \left(\Sigma_h (\bar p, \bar v)\right) (z)}
    \d{z} \,.
  \end{equation}

  Assume that at time $s$ no interaction takes place and choose $h$
  sufficiently small so that in the time interval $[s, s+h]$ no
  interaction takes place and no wave hits any of the lines
  $z = \pm \epsilon^2$, $z = 0$, $z = m\pm \epsilon^2$ and $z=m$.

  We now continue to estimate the right hand side in~\eqref{eq:9}
  limited to $\left]-\infty, 0\right]$. Let $z_1, z_2, \ldots$ be the
  points of jump of the map
  $z \to (p^{\kappa,\epsilon}, v^{\kappa,\epsilon}) (s,z)$.  Denote by
  $\hat\lambda$ an upper bound for the characteristic speeds in the
  gas phase. Then, we have
  \begin{eqnarray}
    \nonumber
    &
    & \int_{-\infty}^0
      \norma{
      (p^{\kappa,\epsilon},v^{\kappa,\epsilon}) (s+h,z)
      -
      \left(\Sigma_h (\bar p, \bar v)\right) (z)}
      \d{z}
    \\
    \label{eq:unouno}
    & =
    & \sum_{z_i < -\epsilon^2} \int_{z_i - \hat\lambda h}^{z_i+\hat\lambda h}
      \norma{
      (p^{\kappa,\epsilon}, v^{\kappa,\epsilon}) (s+h, z)
      -
      \left(
      \Sigma_h \left((p^{\kappa,\epsilon},
      v^{\kappa,\epsilon}) (s)\right)
      \right) (z)
      } \d{z}
    \\
    \label{eq:unodue}
    &
    & +
      \sum_{z_i \in \left]-\epsilon^2,  0\right[} \int_{z_i - \hat\lambda h}^{z_i+\hat\lambda h}
      \norma{
      (p^{\kappa,\epsilon}, v^{\kappa,\epsilon}) (s+h, z)
      -
      \left(\Sigma_h \left((p^{\kappa,\epsilon}, v^{\kappa,\epsilon}) (s)\right)\right) (z)
      } \d{z}
    \\
    \label{eq:unotre}
    &
    & +
      \int_{-\hat\lambda h}^0
      \norma{
      (p^{\kappa,\epsilon}, v^{\kappa,\epsilon}) (s+h, z)
      -
      \Sigma_h \left((\bar p, \bar v) (s)\right)} \d{z} \,.
  \end{eqnarray}
  A standard procedure yields the estimate of~\eqref{eq:unouno} by
  means of~\cite[(ii) in Lemma~9.1]{BressanLectureNotes}, so that
  \begin{displaymath}
    [\eqref{eq:unouno}]
    =
    \O \, \epsilon \, h \, \tv (p^{\kappa,\epsilon} (s); \,]-\infty, -\epsilon^2[\,) \,.
  \end{displaymath}
  Similarly, since all waves in the strip
  $\left]-\epsilon^2, 0\right[$ have speed $\pm 1$, by~\cite[(i) in
  Lemma~9.1]{BressanLectureNotes} we have
  \begin{displaymath}
    [\eqref{eq:unodue}]
    =
    \O \, h \, \tv (p^{\kappa,\epsilon} (s); \,]-\epsilon^2,0[\,) \,.
  \end{displaymath}
  Consider now~\eqref{eq:unotre}. We use~\cite[Point~2) in
  Theorem~2.2]{ColomboGuerra2010} to estimate the difference between
  $(p^{\kappa,\epsilon}, v^{\kappa,\epsilon})$ and
  $\Sigma_h (\bar p, \bar v)$ that are solutions, respectively, to the
  two initial--boundary value problems
  \begin{displaymath}
    \left\{
      \begin{array}{l}
        \partial_t \mathcal{T} (p) - \partial_z v = 0
        \\
        \partial_t v + \partial_z p = 0
        \\
        (p,v) (0,z) = (p^{\kappa,\epsilon}, v^{\kappa,\epsilon}) (s,0)
        \\
        v (t, 0) = v^{\kappa,\epsilon} (s,0)
      \end{array}
    \right.
    \quad \mbox{ and } \quad
    \left\{
      \begin{array}{l}
        \partial_t \mathcal{T} (p) - \partial_z v = 0
        \\
        \partial_t v + \partial_z p = 0
        \\
        (p,v) (0,z) = (p^{\kappa,\epsilon}, v^{\kappa,\epsilon}) (s,0)
        \\
        v (t, 0) = \overline{v}^{\kappa,\epsilon} (s)
      \end{array}
    \right.
  \end{displaymath}
  with the mean value $\overline{v}^{\kappa,\epsilon}$ as defined
  in~\eqref{eq:ukeps}. Then, we
  apply~\cite[Proposition~4.9]{ColomboGuerraIncompressible} to obtain
  \begin{displaymath}
    \begin{array}{rcl@{\qquad}l}
      [\eqref{eq:unotre}]
      & \leq
      & \O \, \hat \lambda \, h \,
        \modulo{
        v^{\kappa,\epsilon} (s,0)
        -
        \overline{v}^{\kappa,\epsilon} (s)}
      &\mbox{by~\cite[Point~2) in Theorem~2.2]{ColomboGuerra2010}}
      \\
      & \leq
      & \O \, \hat \lambda \, h \,
        \tv (v^{\kappa,\epsilon} (s); \mathcal{L})
      & \mbox{by~\eqref{eq:ukeps}}
      \\
      & \leq
      & \O \, h \, \kappa
      & \mbox{by~\cite[Proposition~4.9]{ColomboGuerraIncompressible}.}
    \end{array}
  \end{displaymath}

  Entirely analogous estimates can be applied to bound the similar
  terms on $\left[m, +\infty\right[$. We thus continue~\eqref{eq:9} as
  follows:
  \begin{eqnarray*}
    \norma{U^{\kappa,\epsilon} (s+h)
    - \bar S_h \left(\bar U^{\kappa,\epsilon} (s)\right)}_{\L1 (\reali_+;\reali^4)}
    & \leq
    & \O \, h \,
      \epsilon \,
      \tv (p^{\kappa,\epsilon} (s); \,]-\infty, -\epsilon^2[\, \cup \,]m+\epsilon^2, +\infty [\,)
    \\
    &
    & +
      \O \, h \,
      \tv (p^{\kappa,\epsilon} (s); \,]-\epsilon^2,0[\, \cup \,]m, m+\epsilon^2[\,)
    \\
    &
    & +
      \O \, h \, \kappa \,.
  \end{eqnarray*}

  We pass now to the second addend in~\eqref{eq:questa},
  using~\eqref{eq:10}
  and~\cite[Proposition~4.9]{ColomboGuerraIncompressible},
  \begin{eqnarray}
    \label{eq:qui}
    \!\!\!\!\!\!
    &
    & \norma{
      w^{\kappa,\epsilon} (s+h)
      -
      \left[
      w^{\kappa,\epsilon} (s)
      +
      h F\left(U^{\sigma} ,w^{\kappa,\epsilon} (s)\right)\right]
      }
    \\
    \nonumber
    \!\!\!\!\!\!
    & =
    & \norma{
      \left[
      \begin{array}{c}
        \overline{v}^{\kappa,\epsilon} (s+h)
        -
        \overline{v}^{\kappa,\epsilon} (s)
        -
        h F\left(U^\sigma, w^{k,\epsilon} (s)\right)
        \\
        \overline{v}^{\kappa,\epsilon} (s+h)
        -
        \overline{v}^{\kappa,\epsilon} (s)
        -
        h F\left(U^\sigma, w^{k,\epsilon} (s)\right)
      \end{array}
    \right]
    }
    \\
    \nonumber
    \!\!\!\!\!\!
    & =
    & \sqrt2 \,
      \modulo{
      \overline{v}^{\kappa,\epsilon} (s+h)
      -
      \overline{v}^{\kappa,\epsilon} (s)
      -
      \frac1m \, h \, (p^\sigma_0 - p^\sigma_m)
      }
    \\
    \nonumber
    \!\!\!\!\!\!
    & =
    & \frac{\sqrt2}{m} \,
      \modulo{
      \int_{\mathcal{L}} v^{\kappa,\epsilon} (s+h, z) \d{z}
      -
      \int_{\mathcal{L}} v^{\kappa,\epsilon} (s, z) \d{z}
      -
      \int_s^{s+h}
      \left(
      p^{\kappa,\epsilon} (\sigma,0) - p^{\kappa,\epsilon} (\sigma,m)
      \right) \d\sigma
      }
    \\
    \nonumber
    \!\!\!\!\!\!
    &
    & +
      \frac{\sqrt2 \, h}{m} \,
      \modulo{
      (p^\sigma_0 - p^\sigma_m)
      -
      \left(
      p^{\kappa,\epsilon} (s,0) - p^{\kappa,\epsilon} (s,m)
      \right)
      }
    \\
    \nonumber
    \!\!\!\!\!\!
    & \leq
    & \frac{\sqrt2}{m} \,
      \modulo{
      \int_s^{s+h}
      \frac{\d{~}}{\d{\sigma}}
      \int_{\mathcal{L}} v^{\kappa,\epsilon} (\sigma, z) \d{z} \d\sigma
      -
      \int_s^{s+h}
      \left(
      p^{\kappa,\epsilon} (\sigma,0) - p^{\kappa,\epsilon} (\sigma,m)
      \right) \d\sigma
      }
      + \O \, h \, \kappa
    \\
    \nonumber
    \!\!\!\!\!\!
    & \leq
    & \frac{\sqrt2}{m}
      \Bigl|
      \int_s^{s+h}
      \bigl[
      \sum_{z_i \in \,]0,m[\,}
      \left(
      v^{\kappa,\epsilon} (\sigma, z_i-) - v^{\kappa,\epsilon} (\sigma, z_i+)
      \right)
      \dot z_i
    \\
    \nonumber
    \!\!\!\!\!\!
    &
    & \qquad\qquad\qquad\qquad
      +
      \sum_{z_i \in \,]0,m[\,}
      \left(
      p^{\kappa,\epsilon} (\sigma, z_i+) - p^{\kappa,\epsilon} (\sigma, z_i-)
      \right)
      \bigr] \d{\sigma}
      \Bigr|
      + \O \, h \, \kappa
    \\
    \nonumber
    \!\!\!\!\!\!
    & \leq
    & \frac{\sqrt2}{m} \!
      \int_s^{s+h}
      \!\!\!\!
      \sum_{z_i \in \,]0,m[\,} \!\!
      \modulo{
      \left(
      v^{\kappa,\epsilon} (\sigma, z_i-) - v^{\kappa,\epsilon} (\sigma, z_i+)
      \right)
      \dot z_i
      +
      \left(
      p^{\kappa,\epsilon} (\sigma, z_i+) - p^{\kappa,\epsilon} (\sigma, z_i-)
      \right)
      } \! \d{\sigma}
    \\
    \nonumber
    \!\!\!\!\!\!
    &
    & + \O \, h \, \kappa
  \end{eqnarray}
  We estimate the integral term in the latter term above in different
  ways, depending on the location of $z_i$:
  \begin{eqnarray*}
    &
    & \int_s^{s+h}
      \!\!\!
      \sum_{z_i \in \,]0,\epsilon^2[\,\cup \,]m-\epsilon^2, m[\,}
      \modulo{
      \left(
      v^{\kappa,\epsilon} (\sigma, z_i-) - v^{\kappa,\epsilon} (\sigma, z_i+)
      \right)
      \dot z_i
      +
      \left(
      p^{\kappa,\epsilon} (\sigma, z_i+) - p^{\kappa,\epsilon} (\sigma, z_i-)
      \right)
      } \d{\sigma}
    \\
    & \leq
    & \int_s^{s+h}
      \!\!\!
      \sum_{z_i \in \,]0,\epsilon^2[\,\cup \,]m-\epsilon^2, m[\,}
      \modulo{
      v^{\kappa,\epsilon} (\sigma, z_i-) - v^{\kappa,\epsilon} (\sigma, z_i+)
      }
      +
      \modulo{
      p^{\kappa,\epsilon} (\sigma, z_i+) - p^{\kappa,\epsilon} (\sigma, z_i-)
      }
      \d{\sigma}
    \\
    & =
    & \O \, h \,
      \tv \left(
      p^{\kappa,\epsilon} (s);
      \,]0,\epsilon^2[\,\cup \,]m-\epsilon^2, m[\,\right)
  \end{eqnarray*}
  since in $]0,\epsilon^2[\,\cup \,]m-\epsilon^2, m[$ we have
  $\dot z_i = 1$.  To bound the remaining terms in~\eqref{eq:qui}, let
  $z_i \in \left]\epsilon^2, m-\epsilon^2\right[$,
  use~\cite[Section~4, Lemma~4.1 and
  Formula~(4.3)]{ColomboGuerraIncompressible} and assume that the jump
  at $z_i$ is solved by a $2$-rarefaction:
  \begin{eqnarray*}
    &
    & \modulo{
      -
      \left(
      v^{\kappa,\epsilon} (\sigma, z_i+) - v^{\kappa,\epsilon} (\sigma, z_i-)
      \right)
      \dot z_i
      +
      \left(
      p^{\kappa,\epsilon} (\sigma, z_i+) - p^{\kappa,\epsilon} (\sigma, z_i-)
      \right)
      }
    \\
    & \leq
    & \epsilon \,
      \modulo{
      v^{\kappa,\epsilon} (\sigma, z_i+)
      -
      v^{\kappa,\epsilon} (\sigma, z_i-)}
    \\
    &
    & +
      \Bigg|
      -\frac{1}{\kappa} \,
      \sqrt{-\frac{1}{\mathcal{T}' \left(\Pi_\kappa\left(p^{\kappa,\epsilon} (\sigma, z_i-)\right)\right)}} \,
      \kappa \,
      \left(
      p^{\kappa,\epsilon} (\sigma, z_i+) -p^{\kappa,\epsilon} (\sigma, z_i-)
      \right)
    \\
    &
    & \qquad \qquad \qquad\qquad \qquad \qquad
      \times
      F\left(
      \Pi_\kappa\left(p^{\kappa,\epsilon} (\sigma, z_i-)\right),
      \Pi_\kappa\left(p^{\kappa,\epsilon} (\sigma, z_i+)\right)
      \right)
    \\
    &
    & \qquad
      +
      \left(
      p^{\kappa,\epsilon} (\sigma, z_i+) - p^{\kappa,\epsilon} (\sigma, z_i-)
      \right)
      \Bigg|
    \\
    & \leq
    & \epsilon \, \modulo{v^{\kappa,\epsilon} (\sigma, z_i+) - v^{\kappa,\epsilon} (\sigma, z_i-)}
    \\
    &
    & +
      \frac{
      \modulo{p^{\kappa,\epsilon} (\sigma, z_i+) - p^{\kappa,\epsilon} (\sigma, z_i-)}
      }{
      F\left(
      \Pi_\kappa\left(p^{\kappa,\epsilon} (\sigma, z_i-)\right),
      \Pi_\kappa\left(p^{\kappa,\epsilon} (\sigma, z_i-)\right)
      \right)
      }
    \\
    &
    & \times
      \modulo{
      -F\left(
      \Pi_\kappa\left(p^{\kappa,\epsilon} (\sigma, z_i-)\right),
      \Pi_\kappa\left(p^{\kappa,\epsilon} (\sigma, z_i+)\right)
      \right)
      +
      F\left(
      \Pi_\kappa\left(p^{\kappa,\epsilon} (\sigma, z_i-)\right),
      \Pi_\kappa\left(p^{\kappa,\epsilon} (\sigma, z_i-)\right)
      \right)}
    \\
    & =
    & \epsilon \, \modulo{v^{\kappa,\epsilon} (\sigma, z_i+) - v^{\kappa,\epsilon} (\sigma, z_i-)}
    \\
    &
    & +
      \O
      \modulo{p^{\kappa,\epsilon} (\sigma, z_i+) - p^{\kappa,\epsilon} (\sigma, z_i-)}
      \;
      \modulo{
      \Pi_\kappa\left(p^{\kappa,\epsilon} (\sigma, z_i-)\right)
      -
      \Pi_\kappa\left(p^{\kappa,\epsilon} (\sigma, z_i-)\right)
      }
    \\
    & =
    & \epsilon \, \modulo{v^{\kappa,\epsilon} (\sigma, z_i+) - v^{\kappa,\epsilon} (\sigma, z_i-)}
      +
      \O \,
      \kappa^2 \,
      \modulo{p^{\kappa,\epsilon} (\sigma, z_i+) - p^{\kappa,\epsilon} (\sigma, z_i-)}^2
    \\
    & =
    & \epsilon \, \modulo{v^{\kappa,\epsilon} (\sigma, z_i+) - v^{\kappa,\epsilon} (\sigma, z_i-)}
      +
      \O \,
      \kappa^2 \, \epsilon \,
      \modulo{p^{\kappa,\epsilon} (\sigma, z_i+) - p^{\kappa,\epsilon} (\sigma, z_i-)} \,.
  \end{eqnarray*}
  When dealing with a $2$-shock we obtain the simpler estimate
  \begin{displaymath}
    \modulo{
      -
      \left(
        v^{\kappa,\epsilon} (\sigma, z_i+) - v^{\kappa,\epsilon} (\sigma, z_i-)
      \right)
      \dot z_i
      +
      \left(
        p^{\kappa,\epsilon} (\sigma, z_i+) - p^{\kappa,\epsilon} (\sigma, z_i-)
      \right)
    }
    \leq
    \epsilon \,
    \modulo{v^{\kappa,\epsilon} (\sigma, z_i+) - v^{\kappa,\epsilon} (\sigma, z_i-)}
  \end{displaymath}
  while the cases of waves of the first family are entirely analogous.

  Summarizing:
  \begin{eqnarray*}
    [\eqref{eq:qui}]
    & \leq
    & \O \, h \,
      \tv \left(
      p^{\kappa,\epsilon} (s);
      \,]0,\epsilon^2[\,\cup \,]m-\epsilon^2, m[\,
      \right)
    \\
    &
    & +  \O \, h \, \epsilon \,
      \left(
      \tv\left(v^{\kappa,\epsilon} (s); \,]\epsilon^2, m-\epsilon^2[\, \right)
      +
      \tv\left(p^{\kappa,\epsilon} (s); \,]\epsilon^2, m-\epsilon^2[\, \right)
      \right)
      + \O \, h \, \kappa
  \end{eqnarray*}

  By~\cite[Formula~(4.32) in
  Proposition~4.9]{ColomboGuerraIncompressible} we finally obtain,
  \begin{eqnarray*}
    &
    & d
      \left(
      (U^{\kappa,\epsilon},w^{\kappa,\epsilon}) (s+h),
      \mathcal{F} (h) \left((U^{\kappa,\epsilon},w^{\kappa,\epsilon}) (s)\right)
      \right)
    \\
    & \leq
    & \O \, h \,
      \epsilon \,
      \tv (p^{\kappa,\epsilon} (s); \,]-\infty, -\epsilon^2[\, \cup \,]m+\epsilon^2, +\infty [\,)
    \\
    &
    & +
      \O \, h \,
      \tv (p^{\kappa,\epsilon} (s); \,]-\epsilon^2,0[\, \cup \,]m, m+\epsilon^2[\,)
    \\
    &
    & +
      \O \, h \, \kappa
    \\
    &
    & +
      \O \, h \,
      \tv \left(
      p^{\kappa,\epsilon} (s);
      \,]0,\epsilon^2[\,\cup \,]m-\epsilon^2, m[\,
      \right)
    \\
    &
    & +  \O \, h \, \epsilon \,
      \left(
      \tv\left(v^{\kappa,\epsilon} (s); \,]\epsilon^2, m-\epsilon^2[\, \right)
      +
      \tv\left(p^{\kappa,\epsilon} (s); \,]\epsilon^2, m-\epsilon^2[\, \right)
      \right)
    \\
    & =
    & \O \, h
      \left(
      \epsilon
      +
      \kappa
      +
      \tv \left(p^{\kappa,\epsilon} (s);
      \,[-\epsilon^2, \epsilon^2[\, \cup \,]m-\epsilon^2, m+\epsilon^2[\,\right)
      \right)
  \end{eqnarray*}
  whence, by~\eqref{eq:dX}
  \begin{eqnarray*}
    &
    & d
      \left(
      (U^{\kappa,\epsilon},w^{\kappa,\epsilon}) (t),
      S_t \left((U^{\kappa,\epsilon},w^{\kappa,\epsilon}) (0)\right)
      \right)
    \\
    &
    & \leq
      \O
      \int_0^t
      \left(
      \epsilon
      +
      \kappa
      +
      \tv \left(p^{\kappa,\epsilon} (s);
      \,[-\epsilon^2, \epsilon^2[\, \cup \,]m-\epsilon^2, m+\epsilon^2[\,\right)
      \right)
      \d{s}
  \end{eqnarray*}
  Changing the order of integration and using~\cite[Formula~(4.33) in
  Proposition~4.9]{ColomboGuerraIncompressible}, we get
  \begin{eqnarray*}
    &
    & \int_0^t
      \tv \left(p^{\kappa,\epsilon} (s);
      \,[-\epsilon^2, \epsilon^2[\, \cup \,]m-\epsilon^2, m+\epsilon^2[\,\right)
      \d{s}
    \\
    & =
    & \int_{\,[-\epsilon^2, \epsilon^2[\, \cup \,]m-\epsilon^2, m+\epsilon^2[\,}
      \tv\left(p^{\kappa,\epsilon} (\cdot, z); [0,t]\right) \d{z}
    \\
    & =
    & \O \, \frac{\epsilon^2}{\kappa}
  \end{eqnarray*}
  so that
  \begin{displaymath}
    d
    \left(
      (U^{\kappa,\epsilon},w^{\kappa,\epsilon}) (t),
      S_t \left((U^{\kappa,\epsilon},w^{\kappa,\epsilon}) (0)\right)
    \right)
    =
    \O \left((\epsilon+\kappa)t + \frac{\epsilon^2}{\kappa}\right)
  \end{displaymath}
  Using~\eqref{eq:agognata}, the proof of~1.~is completed.

  Hence, the trajectory of the semigroup $S$ with initial datum
  $(U,w) (0)$ is defined for all $t \in \reali_+$.
\end{proofof}

\noindent\textbf{Acknowledgment:} The present work was supported by
the PRIN~2012 project \emph{Nonlinear Hyperbolic Partial Differential
  Equations, Dispersive and Transport Equations: Theoretical and
  Applicative Aspects} and by the GNAMPA~2014 project
\emph{Conservation Laws in the Modeling of Collective Phenomena}.

\end{document}